\documentclass[a4paper, 11pt]{amsart}
\usepackage{amsfonts, amsmath, amssymb, graphics, epsfig}
\usepackage{xcolor}
\usepackage{enumerate, amstext, amsthm,bm}

\theoremstyle{plain}
\newtheorem{theorem}{Theorem}[section]
\newtheorem{lemma}[theorem]{Lemma}
\newtheorem{corollary}[theorem]{Corollary}
\newtheorem{proposition}[theorem]{Proposition}

\theoremstyle{definition}

\newtheorem{example}[theorem]{Example}
\theoremstyle{remark}

\def\IC{{\mathbb C}}

\def\IR{{\mathbb R}}

\def\b0{{\mathbf 0}}

\def\cB{{\mathcal B}}

\def\cF{{\mathcal F}}

\def\cK{{\mathcal K}}

\def\cS{{\mathcal S}}

\def\cH{{\mathcal H}}

\def\BH{{{\mathcal B}(\cH)}}
\def\BK{{{\mathcal B}(\cK)}}

\def\tr{{\rm tr}\,}

\def\rank{{\rm rank}\,}
\def\diag{{\rm diag}\,}

\def\[{\left [}
\def\]{\right ]}
\def\({\left (}
\def\){\right )}

\def\la{{\langle}}
\def\ra{{\rangle}}

\def\<{{\langle}}
\def\>{{\rangle}}
\def\1{{\mathbf 1}}

\def\la{{\langle}}
\def\ra{{\rangle}}

\newcommand{\vertiii}[1]{{\left\vert\kern-0.25ex\left\vert\kern-0.25ex\left\vert #1
    \right\vert\kern-0.25ex\right\vert\kern-0.25ex\right\vert}}

\begin{document}
\openup .85\jot
\title[Unitarily Invariant Norms on Operators]
{Unitarily invariant Norms on Operators}

\author[Jor-Ting Chan]{Jor-Ting Chan}
\author[Chi-Kwong Li]{Chi-Kwong Li}
\address[1]{Email: jtchan@hku.hk}
\address[2]{Department of Mathematics, The College of William
\& Mary, Williamsburg, VA 13185, USA. Email: ckli@math.wm.edu}

\begin{abstract}
Let $f$ be a symmetric norm on $\IR^n$ and let $\BH$ be the set of all bounded linear operators
on a Hilbert space $\cH$ of dimension at least $n$.
Define a norm
on $\BH$ by
$\|A\|_f = f(s_1(A), \dots, s_n(A))$, where
$s_k(A) = \inf\{\|A-X\|: X\in \BH \hbox{ has rank less than } k\}$
is the $k$th singular value of $A$.
Basic properties of the norm $\|\cdot\|_f$
are obtained including
some norm inequalities and characterization of the
equality case. Geometric properties
of the unit ball of the norm
are obtained; the results are used
to determine the structure of maps $L$ satisfying
$\|L(A)-L(B)\|_f=\|A - B\|_f$ for any $A, B \in \BH$.
\end{abstract}

\date{}

\maketitle

\section{Introduction}
Let $\mathcal H$ be a complex Hilbert space
and $\BH$ the algebra of all bounded linear operators
on $\mathcal H$.
If $\cH$ has finite dimension
$n$, we can identify
$\cH$ and $\BH$ respectively as $\IC^n$ and
$M_n$, the set of all $n\times n$ complex matrices.
Let $\|\cdot\|_{\rm sp}$ be the operator norm on $\BH$. For any $A\in\BH$, let
$$s_k(A) = \inf\{\|A-X\|_{\rm sp}: X\in \BH \hbox{ has rank less than } k\}$$
be the $k$th singular value of $A$. They are also called the $s$-numbers or approximation numbers
of $A$.  See
\cite{P} for the background.
If $A$ is compact, then $s_k(A)$ is the $k$th largest eigenvalue of $(A^*A)^{1/2}$.
We can use the singular values of matrices to define
those for a general operator $A\in \BH$ as
$$s_k(A) = \sup \{s_k(X^*AY): X, Y: \mathbb C^n\to \cH, X^*X = Y^*Y = I_k\}.$$

\medskip
A norm $f$ on $\IR^n$ is called symmetric if
$f(Px) = f(x)$ for any $x \in \IR^n$ and any permutation or
diagonal orthogonal matrix $P$.
Suppose $f$ is a symmetric norm on $\IR^n$ and
$n \le \dim \cH$.
One can define
a norm on $\BH$ by
$$\|A\|_f = f(s_1(A), \dots, s_n(A))\quad\mbox{for all $A$}.$$
For example,
if $f(x) = \sum_{j=1}^n (|c_j x_j|^p)^{1/p}$ for
$x = (x_1, \dots, x_n)^t\in \IR^n$, where
$(c_1, \dots, c_n)$ is a nonzero vector, then $\|\cdot\|_f$ is
the $(c,p)$-norm studied in \cite{CLT}.

\medskip
It is clear that $\|\cdot\|_f$ is unitarily invariant, i.e.,
$$\|UAV\| = \|A\| \quad \hbox{ for all } A \hbox{ and unitary }
U, V \hbox{ in } \BH.$$
Unitarily invariant norms on $M_n$ are well studied.
There has also been considerable interest in unitarily invariant
norms on $\BH$, but most of the study were
restricted  to symmetrically norm ideals in $\BH$
to ensure that the norm
is bounded; for example, see \cite{A, FL, GK, S, So}.

\medskip 
In this paper, we study norm properties and distance preserving maps related to 
$\|\cdot\|_f$ defined above. 
In Section 2, we study basic properties for
$\|\cdot\|_f$, and obtain results analogous to those 
on matrices and compact operators; e.g., see
\cite{DLL, JL, LT, M}.
In Section 3, we obtain some geometrical properties of the unit
norm ball of $\|\cdot\|_f$ and use the result to determine
the structure of maps $L:\BH\rightarrow \BH$ such that
$\|A-B\|_{f} = \|L(A)-L(B)\|_{f}$ for all $A, B \in \BH$.
These extend the results on isometries and  distance preserving maps
of unitarily invariant norms in the literature;
e.g., see \cite{A, CLS, CLT, GM, Ka, Li2, LT2, So}.

\section{Basic norm properties and inequalities}

The construction of $\|\cdot\|_f$ on $\BH$ also arises naturally from other
considerations.
For any norm $\|\cdot\|$ on $M_n$ with
$n \le \dim\cH$, one can define a
norm on $\BH$ by
$$\|A\|_{\rm ui1} = \sup\{ \|X^*AY\|:  X, Y: \mathbb C^n\to \cH, X^*X = Y^*Y = I_n\}.$$
Then $\|\cdot\|_{\rm ui1}$ is unitarily invariant.
Clearly, if $n = 1$
and $\|[a]\| = |a|$, then $\|\cdot\|_{\rm ui1}$ is the operator norm.
In fact, one may assume that the norm $\|\cdot\|$ on $M_n$ is
unitarily invariant. Else, we may modify the norm
as
$$\vertiii A = \max\{\|UAV\|: U, V \in M_n \hbox{ are unitary}\}.$$
Then $\vertiii{\ \cdot\ }$ is unitarily invariant;  $\|\cdot\|$ and
$\vertiii{\ \cdot\ }$ induce the same
$\|\cdot\|_{\rm ui1}$ on $\BH$.
By a result of von Neumann in
\cite{vN} (see  also \cite{LT2, M}),

\medskip
(U1) there is a symmetric norm $f$ on $\IR^n$
such that $\|A\|_{\rm ui1} = f(s_1(A), \dots, s_n(A))$ for all $A\in M_n$.

\medskip\noindent
As a result, $\|A\|_{\rm ui1} = \|A\|_f$ for all $A \in \BH$.
So, instead of using a symmetric norm  $f: \IR^n \rightarrow \IR$ to define
$\|\cdot\|_f$ on $\BH$, one can start with an arbitrary norm $\|\cdot\|$ on
$M_n$ and the corresponding $\|\cdot\|_{\rm ui1}$ defined above will coincide with
$\|\cdot\|_f$ for some symmetric norm $f$.

\medskip
Also, let $\|\cdot\|$ be any norm on $\BH$.
If $n \le \dim \cH$,
one can define
$$\|A\|_{\rm ui2} = \sup\{ \|X^*AY\|: X, Y\hbox{ are rank $n$ partial isometries}\}.$$
Then $\|\cdot\|_{\rm ui2}$ is a unitarily invariant norm
on $\BH$ such that
$$
\|A\| \le \|A\|_{\rm ui2} \qquad \hbox{ for all } A \in \BH
\hbox{ with } \rank(A) \le n.
$$
Take two orthonormal families $\{x_1, \dots, x_n\}$ and $\{y_1, \dots, y_n\}$
in $\cH$ and define $f: \IR^n \rightarrow \IR$ by
$$f(a_1, \dots, a_n) = \left\|\mbox{$\sum_{j=1}^n a_j x_j y_j^*$}\right\|_{\rm ui2}\quad\mbox{for all $(a_1, \ldots, a_n)\in \IR^n$}.$$
Then one can check that $f$ is a symmetric norm on $\IR^n$ and
$$\|A\|_{\rm ui2} = \|A\|_f \qquad \hbox{ for  all } A  \in \BH.$$
So, $\|\cdot\|_f$ can be generated from any norm on $\BH$.

There are other interesting connections between unitarily invariant norms and symmetric
norms in the finite dimensional case,
see \cite{Li2} and \cite{LT}.
In particular,
if $\|\cdot\|$ is a unitarily
invariant norm on $M_n$, then (U1) holds, and

\medskip
\begin{itemize}
\item[(U2)] there is a compact set $\cS \subseteq \IR^n$
consisting of
vectors
$c = (c_1, \dots, c_n)$ with $c_1 \ge \cdots \ge c_n \ge 0$
such that
$$\|A\| = \max\{ \|A\|_c: c \in \cS\},$$
where $\|A\|_c = \sum_{j=1}^n c_j s_j(A)$ is the $c$-norm of $A$.
\end{itemize}

\medskip
From the above discussion we have the following proposition, which
will be useful in our study.

\begin{proposition} \label{uinorm}
Let $f: \IR^n\rightarrow \IR$ be a symmetric norm.
Suppose $\cH$ is a Hilbert space with $\dim \cH \ge n$,
and  $\|\cdot\|$  is a norm on $\BH$ defined by
\begin{equation}\label{ui-norm}
\|A\| = f(s_1(A), \dots, s_n(A)) \quad \hbox{ for all } A \in \BH.
\end{equation}
Then there is a compact set $\cS \subseteq \IR^n$
consisting of
vectors
$c = (c_1, \dots, c_n)$ with $c_1 \ge \cdots \ge c_n \ge 0$
such that
\begin{equation}
\label{c-norm}
\|A\|
= \max\{ \|A\|_c: c \in \cS\}.
\end{equation}
\end{proposition}

\medskip
By the above result,
we can use the properties of symmetric norms and
$c$-norms to study  unitarily invariant norm $\|\cdot\|$ defined
on $\BH$ as in  (\ref{ui-norm}).
For example, it is easy to extend \cite[Theorem 1]{JL}
concerning the comparison of a unitarily invariant norm
and the numerical radius defined by
$$r(A) = \sup\{ |\la Ax, x\ra|:
x\in \cH, \la x, x\ra = 1\}.$$
To present our result and proof, we need the following definition.
For two nonnegative vectors $x, y \in \IR^n$
we say that  $x$ is weakly majorized by
$y$, denoted by $x \prec_w y$,
if the sum of the $k$ largest entries of $x$ is not
larger than that of $y$ for every $1\le k\le n$. See \cite{MO}.

\begin{proposition}
Suppose $\|\cdot\|$ is a norm
on $\BH$ defined as in {\rm (\ref{ui-norm})}.
Let
$$\alpha = f(e_1) \quad \hbox{ and }
\quad \beta =
\begin{cases}
2f(e_1 + \cdots + e_k) & \hbox{ if } n = 2k,\cr
2f(e_1 + \cdots + e_k+e_{k+1}/2) & \hbox { if } n = 2k+1.\cr
\end{cases}$$
Then for any $A \in \BH$,
$$\alpha r(A) \le \|A\| \le \beta r(A).$$
The constants $\alpha$ and $\beta$ are the best possible in view of the following:
\begin{itemize}
\item
$\alpha r(A) = \alpha s_1(A) = \|A\|$ if $A = xx^*$ for a unit vector $x$;
\item
$\|A\| = \beta r(A)$ if
$A = \begin{cases} \sum_{j=1}^k 2x_j y_j^* & \hbox{ if } n = 2k,\cr
\sum_{j=1}^k 2x_j y_j^* + x_{k+1}x_{k+1}^* & \hbox{ if } n = 2k+1,\cr
\end{cases}$ \\
for an orthonormal set
$\cF = \{x_1,\dots, x_k, y_1, \dots, y_k\} \subseteq \cH$ and $x_{k+1}\in
\cF^\perp$ if $n = 2k+1$.
\end{itemize}
\end{proposition}

\it Proof. \rm
Denote by $\{E_{jk}\}_{j, k = 1}^n$
the standard basis for $M_n$.
Let $s(Y)$ be the vector of singular values of
$Y \in M_n$.
If $A \in \BH$ has $r(A) = 1$, then for any
$X, Y: \IC^n\to \cH$ such that $X^*X = Y^*Y = I_n$, we have
from the proof of Theorem 1 in \cite{JL},
$$s(E_{11}) \prec_w s(X^*AY)
\prec_w s(F),$$
where $F = 2(E_{12} + E_{34} + \cdots + E_{2k-1,2k})$ if $n = 2k$
and $F = 2(E_{12} + E_{34} + \cdots  + E_{2k-1,2k}) + E_{nn}$
if  $n = 2k+1$.
It is known that if $R,T \in M_n$ with $s(R) \prec_w s(T)$
then
$|R| \le |T|$ for any unitarily invariant norm $|\cdot|$;
see \cite{LT, M}. We get the inequalities.
The equality cases can be checked directly. \qed

\medskip
We extend the definitions in (\cite[pp. 54-55]{Li2})
for $\|\cdot\|$ defined as in (\ref{ui-norm}):

\medskip

\
(i) $\|\cdot\|$ is a uniform norm if
$\|AB\| \le s_1(A)\|B\|$ and $\|AB\|
\le s_1(B) \|A\|$,

\ (ii) $\|\cdot\|$ is a cross norm if
$\|A\| = s_1(A)$ for any rank one operator
$A \in \BH$.

\ (iii) $\|\cdot\|$ on $\BH$ is submultiplicative if
$\|AB\| \le \|A\|\|B\|$
for all $A, B\in \BH$;

\qquad  if, in addition,
$\|I\|= 1$, it is an algebra (or a ring) norm.

\noindent
We have the following results.

\begin{proposition} \label{2.2}
Suppose $\|\cdot\|$ is a norm
on $\BH$ defined as in {\rm (\ref{ui-norm})}.
\begin{itemize}
\item[(a)] Then $\|\cdot\|$ is a uniform norm.
\item[(b)] The norm $\|\cdot\|$ is a cross norm
if and only if $f(e_1) = 1$.
\end{itemize}
\end{proposition}

\it Proof. \rm (a) As observed above, we may assume that $\|\cdot\|$ is induced by
a unitarily invariant norm $\vertiii{\ \cdot\ }$ on $M_n$. For any $A, B\in \BH$ and $X, Y: \IC^n\to \cH$ such that
$X^*X = Y^*Y = I_n$, let $P$ be any rank $n$ partial isometry whose image contains that of $BY$.
Then by \cite[Lemma 2.2]{DLL},
$$\vertiii {X^*ABY} = \vertiii{(X^*AP)(P^*BY)}\le s_1(X^*AP)\vertiii{P^*BY}\le s_1(A)\|B\|.$$
The result follows.

(b) Evidently, $f(e_1) =1$
if and only if $1 =\|xy^*\| = f(e_1)$
for any unit vectors $x,y$. \qed

\begin{proposition} \label{2.2b}
Suppose $\|\cdot\|$ is a norm
on $\BH$ defined as in {\rm (\ref{ui-norm})}.
Then  $\|\cdot\|$ is submultiplicative if and only if
any one of the following holds.
\begin{itemize}
\item[{\rm (1)}] $\|A\| \ge s_1(A)$ for all $A \in \BH$.
\item[
{\rm (2)}] There are unit vectors $x,y \in \cH$ such that
$\|xy^*\| \ge 1$.
 \item[
{\rm (3)}] $f(e_1) \ge 1$, equivalently,
there is $(c_1, \dots, c_n) \in \cS$ with $c_1 \ge 1$.
\end{itemize}
Moreover,  $\|\cdot\|$ is an algebra norm if and only
if it is the operator norm.
\end{proposition}

\it Proof. \rm
The equivalence of (1) -- (3) can be proved as in the finite
dimensional case, e.g., see \cite{LT}.
Suppose $\|\cdot\|$ is multiplicative. Then
for any unit vector $x \in \cH$,
$\|xx^*\|^2 \ge \|(xx^*)^2\| = \|xx^*\|$
so that $1\le \|xx^*\| = f(e_1)$.  Thus, (3) holds.
Suppose (1) holds. Then by Proposition 2.3 (a), $\|AB\|
\le s_1(A)\|B\| \le \|A\| \|B\|$ for any
$A, B \in \BH$.

To prove the last assertion,
suppose $\|\cdot\|$ is an algebra norm.
Then for any $c = (c_1, \dots, c_n) \in \cS$,
$1 =\|I\| \ge \|I\|_c = \sum_{j=1}^n c_j$. For any $A\in \BH$,
there is a $c\in \cS$ such that
$$\|A\| = \|A\|_c = \sum_{j=1}^n c_js_j(A)\le s_1(A).$$
But $\|\cdot\|$ is submultiplicative, so that $\|A\| \ge s_1(A)$.
Thus, $\|\cdot\|$ is the operator norm. The converse is clear.
\qed

One may obtain similar inequalities for the triple product on $\BH$
defined by
$$A \circ B \circ C = \frac{1}{2}(AB^*C+CB^*A).$$

\begin{proposition} Suppose $\|\cdot\|$ is a norm
on $\BH$ defined as in {\rm (\ref{ui-norm})}.
The following condition are equivalent.
\begin{itemize}
\item[{\rm (1)}]
$\|A\circ B \circ C\| \le \|A\| \|B\| \|C\|$ for all $A, B, C \in
\BH$.

\item[{\rm (2)}] $\|A\circ A \circ A\| \le \|A\|^3$ for all $A, \in
\BH$.

\item[{\rm (3)}] $f(e_1) \ge 1$, equivalently,
there is $(c_1, \dots, c_n) \in \cS$ with
$c_1 \ge 1$.
\end{itemize}
\end{proposition}

\it Proof. \rm If (1) holds, then (2) clearly holds.
If (2) holds, consider $A = xy^*$ for unit vectors $x,y$.
Clearly, $\|A\| = f(e_1)$ and $A\circ A \circ A= A$.
Thus, the inequality in (2) becomes $f(e_1)\le f(e_1)^3$, and hence
$f(e_1)  \ge 1$.

\medskip Let $(s_1,
\dots, s_n), (a_1, \dots, a_n), (b_1, \dots, b_n),
(c_1, \dots, c_n)$ be the vectors of the
$n$ largest singular values of
$A\circ B \circ C, A, B, C,$ respectively.
Then
$\sum_{j=1}^\ell s_j \le \sum_{j=1}^\ell
a_jb_jc_j$ for $\ell = 1, \dots, n$. This ensures that
$$f(s_1, \dots, s_n)\le f(a_1b_1c_1, \dots, a_nb_nc_n).$$
If (3) holds, then for vectors $(x_1, \dots, x_n),
(y_1, \dots, y_n)$ with nonnegative entries
arranging in descending order, we have
$$f(x_1y_1, \dots, x_ny_n) \le f(x_1, \dots, x_n) f(y_1, \dots, y_n).
$$
Consequently,
$$\|A\circ B \circ C\|_f
= f(s_1, \dots, s_n) \le f(a_1b_1c_1, \dots, a_nb_nc_n)
$$
$$
\le f(a_1,\dots, a_n) f(b_1, \dots, b_n) f(c_1, \dots, c_n)
= \|A\|_f \|B\|_f \|C\|_f.$$
Thus (1) holds. \qed

\medskip
The next theorem extends the result for the finite dimensional
case in  \cite{DLL}. In fact, there was a gap in line 4 of
the proof of  \cite[Theorem 2.4]{DLL}, which
asserted that $1 \le s_1(AB)/\|A\|\|B\|$
without justification.
Our proof fills the gap.

\begin{theorem} \label{2.6}
Suppose $\|\cdot\|$
on $\BH$ defined as in {\rm (\ref{ui-norm})}
is a submultiplicative norm.
Then there exist nonzero $A, B \in \BH$
such that $\|AB\| = \|A\|\|B\|$ if and only if
$\|\cdot\|$ is a cross norm.

\noindent
In case $\|\cdot\|$ is a cross norm, the following conditions
hold.
\begin{itemize}
\item[
(1)] Two nonzero operators $A, B \in \BH$ satisfy
$\|AB\|=  \|A\| \|B\|$ if and only if
$s_1(AB) = \|A\|\|B\|$.
\item[
(2)]
Assume that
$\{T\in \BH: \|T\| = s_1(T)\}$
is a subset of rank one matrices.
Two nonzero operators $A, B \in \BH$ satisfy
$\|AB\|=  \|A\| \|B\|$ if and only if
$A = xy^*$ and $B = y^*z$ for some nonzero
$x,y,z \in \cH$ with $\la y, y\ra = 1$.
\end{itemize}
\end{theorem}

\it Proof. \rm
Suppose $\|\cdot\|$ is submultiplicative.
If $\|\cdot\|$ is a cross norm, then for any unit vector
$x \in \cH$ and $A = B= xx^*$, we have
$\|AB\| = \|A\| \|B\| = 1.$

Suppose $\|\cdot\|$ is not a cross norm. Then
$\|xy^*\| = d > 1$ for any unit vectors $x, y \in \cH$.
We will show that $\|A\|\|B\| > \|AB\|$
for any nonzero $A, B \in \BH$. For any $\varepsilon\in (0, s_1(A))$,
take a unit vector $x\in \cH$ such that $\|Ax\| > s_1(A) - \varepsilon$
and let $y = Ax/\|Ax\|$. Then $s_1(Axy^*) \le s_1(A)$
and $s_j(Axy^*) = 0$ for all $j > 1$. Consequently,
$\sum_{j = 1}^k s_j(A)\ge \sum_{j = 1}^k s_j(Axy^*)$ for all $k = 1, \ldots, n$.
It follows that $\|A\| \ge \|Axy^*\| = \|Ax\| \|yy^*\| > (s_1(A) - \varepsilon)d$.
As this is true for all such $\varepsilon$, $\|A\|\ge s_1(A)d > s_1(A)$.
Thus $\|A\| \|B\| > s_1(A) \|B\| \ge \|AB\|$.

\medskip
Now, assume that $\|\cdot\|$ is a cross norm.

(1) Note that
$s_1(AB)  \le \|AB\| \le \|A\| \|B\|$.
If $s_1(AB) = \|A\| \|B\|$, then $\|AB\| =\|A\| \|B\|$.

Conversely, if $\|AB\| =\|A\| \|B\|$,  then
by Proposition 2.3 (a) and Proposition 2.4 (1), $s_1(A) =\|A\|$ and $s_1(B) = \|B\|$.
Suppose $s_1(AB) < \|A\|\|B\| = s_1(A)s_1(B)$.
Also, recall that $s_j(AB) \le s_j(A)s_1(B)$
for $j = 2, \dots, n$.
Hence, for $k  = 1, \dots, n$,
$$\sum_{j=1}^k s_j(AB) < \sum_{j=1}^k s_j(A)s_1(B)
\quad \hbox{ so that } \quad
\sum_{j=1}^k s_j(AB) \le t \sum_{j=1}^k s_j(A)s_1(B) = \sum_{j=1}^k s_j(ts_1(B)A)$$
for some  $t \in (0,1)$.
Consequently,
$$\|AB\| \le \|t s_1(B) A\|
= t \|A\| \|B\| < \|A\| \|B\|,$$
 which
is a contradiction.

(2) Direct verification.  \qed

\medskip
We remark that in Theorem \ref{2.6} (2), the assumption that
$\{T\in \BH: \|T\| = s_1(T)\}$
is a subset of rank one matrices is equivalent to the set equality
$$\{T\in \BH: \|T\| = s_1(T)=1\} = \{xy^*: x, y \in \cH \hbox{ are unit vecotrs}\}.$$
By \cite[Example 2.5]{DLL}, the result in Theorem \ref{2.6}
(2) is not true if the above set equality fails.
For instance, if $\|\cdot\|$ is the spectral norm, any unitary $A, B$ will
satisfy $\|AB\| = \|A\| \|B\|$.

\medskip
Specializing the above results to the $c$-norm, we have
the following.

\begin{corollary}
Suppose $c = (c_1, \dots, c_n)$ has positive entries arranged
in descending order.
\begin{itemize}
\item[{\rm (a)}]
The norm $\|\cdot\|_c$
is submultiplicative if and  only if $c_1 \ge 1$.
\item[{\rm (b)}]
The norm $\|\cdot\|_c$ is a cross norm if and only if
$c_1 = 1$.
\item[{\rm (c)}] Suppose $c_1\ge 1$.
There exist nonzero $A, B \in \BH$ such that
$\|AB\|_c = \|A\|_c\|B\|_c$
if and only if $c_1 = 1$.
\item[{\rm (d)}] Suppose $c=(1,0,\dots,0)$.
Two nonzero operators $A, B \in \BH$
satisfy $\|AB\|_c = \|A\|_c\|B\|_c$
if and only only if $s_1(AB) = s_1(A)s_1(B)$.
\item[{\rm (e)}] Suppose $c_1 = 1$ and
$c_2 > 0$. Two nonzero operators $A, B \in \BH$ satisfy
$\|AB\|_c = \|A\|_c \|B\|_c$ if and only if
$A = xy^*$ and $B = yz^*$ for some nonzero vectors
$x, y, z \in \cH$.
\end{itemize}
\end{corollary}

\section{Distance preserving maps}

In this section, we always assume that
$\|\cdot\|$ is defined as in
(\ref{ui-norm}) for $n > 1$. So,
$$
\|A\| = f(s_1(A), \dots, s_n(A))\ (= \|A\|_f)
= \max\{\|A\|_c: c \in \cS\},
$$
for a symmetric norm $f$ on $\IR^n$ and a compact set $\cS$ in $\IR^n$.
This property will be used later without further reference.
Our main result in this section is the following.

\begin{theorem} \label{isometry} A surjective map
$L:\BH \rightarrow \BH$ satisfies
\begin{equation}\label{normpreserving}
\|A-B\| = \|L(A)-L(B)\| \quad \hbox{ for all } A, B
\in  \BH
\end{equation}
if and only if there are unitary operators $U, V \in \BH$ such that
$L$ has the form
$A \mapsto U^*\phi(A)V + R_0$, where $R_0\in \BH$ and
$\phi$ is one of the maps
$$A \mapsto A,  \quad
A \mapsto A^t,  \quad
A \mapsto A^*,  \quad
A \mapsto (A^*)^t.$$
Here $A^t\in \BH$ is the transpose of $A$ with respect to a fixed
orthonormal basis of $\cH$.
\end{theorem}

If $\cH$ is finite dimensional, then the result
holds without the assumption that $L$ is surjective.
It is a special case of \cite[Theorem 3.2]{CLS}.
We will therefore focus on the case when
$\cH$ has infinite dimension.
Note also that a characterization of complex isometric isomorphisms
on $\BH$, when $\dim \cH$ is infinite, is given in \cite{CLT}
for a special class of unitarily invariant norms.

\medskip
Let $L:\BH \rightarrow \BH$ be a surjective map satisfying the hypothesis of Theorem \ref{isometry}.
Then by the Mazur-Ulam Theorem, the map $\hat L: A \mapsto L(A)-L(0)$ is
real linear, surjective, and satisfies
$\|\hat L(A)\| = \|A\|$ for all $A \in \BH$.
So, $\hat L$ preserves $\cB = \{X \in \BH: \|X\| \le 1\}$,
the closed unit ball for $\|\cdot\|$, and the extreme points of it.
In the following, we give some auxiliary results concerning the
extreme points of $\cB$. In particular, we show that
$U/\|U\|$ is an extreme point of $\cB$ for
any maximal partial isometry $U$, i.e., $U\in\BH$ such that
$UU^* = I$ or $U^*U = I$.

\begin{lemma} \label{Extreme}
Let $A \in \BH$ be an extreme point
of $\cB$. Then
$A = \sum_{j = 1}^n s_j(A) x_jy_j^* + s_n(A)U$, where
$s_1\ge \cdots\ge s_n\ge 0$,
$\{x_1, \ldots, x_n\}$, $\{y_1, \ldots, y_n\}$ are orthonormal sets in $\cH$
and $U$ is a partial isometry from $\{y_1, \ldots, y_n\}^\perp$ onto $\{x_1, \ldots, x_n\}^\perp$.

\medskip
On the other hand, $U/\|U\|$ is an extreme point of $\cB$
for every maximal partial isometry $U\in \BH$.
\end{lemma}

\it Proof. \rm Let $A \in \BH$ be an extreme point of $\cB$.
We shall show that $A$ is of the given form using the ideas in \cite{CLT}.
First of all, we have $s_k(A) = s_n(A)$ for every $k > n$.
Otherwise, there is a $k > n$ such that $s_k(A) < s_n(A)$
and $s_k(A)$ is an eigenvalue of $|A|$. Take a unit eigenvector $x$
of $|A|$ corresponding to $s_k (A)$ and let $B = (Ax)x^*$. For sufficiently small
$\varepsilon > 0$, $A\pm \varepsilon B$ have the same $n$ largest singular values as $A$.
So, $\|A\pm \varepsilon B\| = 1$, contradicting the choice of $A$.
Next, we claim that $\sigma(|A|)\subseteq \{0, s_1(A), \ldots, s_n(A)\}$.
An argument as above shows that it is true if every non-zero element of $\sigma(|A|)$ is an eigenvalue.
Otherwise we can use the construction in \cite[Lemma 3.7]{CLT} to get a contradiction
if the claim is not true.
So, $A$ is of the form $\sum_{j = 1}^n s_j(A) x_jy_j^* + s_n(A)U$
for orthonormal sets $\{x_1, \ldots, x_n\}$, $\{y_1, \ldots, y_n\}$
and a partial isometry $U$. This $U$ must be a maximal partial isometry.
Otherwise, a contradiction can be obtained as above by considering $B = yx^*$
for unit vectors $x\in {\rm Ker}\,U$ and $y\in {\rm Im}\,U$.

\medskip
The last statement is a consequence of \cite[Lemma 3.5]{CLT}
which asserts that if $U$ is a maximal partial isometry, then
$U/\|U\|_c$ is an extreme point of $\cB_c$, the closed unit ball of
$\|\cdot\|_c$, for every $c\in\cS$. So, take
$c \in \cS$ such that $\|U\| = \|U\|_c$.
Then $U/\|U\| = U/\|U\|_c$ is an extreme point of $\cB_c$.
It must be an extreme point of $\cB$ as $\cB\subseteq \cB_c$.
The lemma is proved.
\qed

\medskip
To prove Lemma \ref{Claim}, we need the following result.

\begin{lemma} \label{Assertion}
Suppose $c = (c_1, \dots, c_p)$ has all entries positive and
$C = \diag(c_1, \ldots, c_p)$. If $R \in M_p$ satisfies
$\|R\|_c = \tr(CR),$
then $R$ is positive semidefinite.
\end{lemma}

\it Proof. \rm  It is known that $R \in M_p$
$R$ is positive semidefinite if and only if
$\tr R = \sum_{j = 1}^p s_j(R)$; e.g., see
\cite[Corollary 3.2]{Li}. In the following, we will show that
$\tr R = \sum_{j = 1}^p s_j(R)$ if $\|R\|_c = \tr(CR)$.

As some of the entries in $c$ may be equal, write $C = \xi_1 I_{p_1} \oplus \cdots
\oplus \xi_{u} I_{p_u}$ with $\xi_1 > \cdots > \xi_u > 0$
and $R = (R_{ij}) \in M_p$ with $R_{jj}
\in M_{p_j}$ for $j = 1, \dots, u$.
It is known if $R$ has diagonal entries
$d_1, \dots, d_p$, then
$|\sum_{j=1}^\ell d_j| \le \sum_{j=1}^{\ell} s_j(R)$;
e.g., see \cite{Li}.
So,
$$\tr(CR)
= \sum_{j=1}^u \xi_j \tr R_{jj}
= \sum_{j=1}^{u-1} (\xi_j-\xi_{j+1}) \tr(R_1 + \cdots + R_j)
+ \xi_u\tr R,$$
and
$$\left|\sum_{j=1}^{u-1} (\xi_j-\xi_{j+1}) \tr(R_1 + \cdots + R_j)
+ \xi_u\tr R\right| \hskip .5in \ $$
$$\ \hskip .3in \le
\sum_{j=1}^{u-1} (\xi_j-\xi_{j+1})
[s_1(R) + \cdots + s_{p_1+ \cdots + p_j}(R)]
+ \xi_u \sum_{j=1}^p s_j(R) = \|R\|_c.$$
Thus, $\tr (CR) = \|R\|_c$  ensures
$\tr R = \sum_{j = 1}^p s_j(R)$, which is
the desired equality mentioned at the beginning of the proof.
Hence $R$ is positive semidefinite.
\qed

\begin{lemma} \label{Claim}
Let $A = a_1 I_{r_1} \oplus \cdots \oplus a_{m}I_{r_m} + a_{m+1}U \in \BH$
with $a_1 > \cdots > a_{m+1} \ge 0$ and $U$ a maximal partial isometry from
$[\ker(A-a_1I)\oplus \cdots \oplus \ker(A-a_mI)]^\perp$ into itself.
Suppose $A = B + D$ and they satisfy
\begin{equation}
\label{norm}
\|\mu B + \nu D\| = \max\{|\mu|, |\nu|\}
\quad\mbox{for any\quad $\mu, \nu \in \IR$.}\end{equation}
Then there is a vector $c\in \cS$ such that $\|A + B\|_c = \|A + B\| = 2$.
Let $B_1$ be the compression of
$B$ onto $\cK = \ker(A-a_1I)$. If $c$ has $p$ nonzero entries, then
one of the following holds.
\begin{itemize}
\item[{\rm (a)}] $p < r_1$ and there is a linear functional $\varphi\in \BK^*$, the dual space of $\BK$, such that
$\varphi(B_1) = 1$. Moreover, if
$Y\in \BH$ satisfies
$\|\mu B + \nu Y\| = \max\{|\mu|, |\nu|\}$
for any $\mu, \nu \in \IR$ and $Y_1$ is the compression of $Y$ onto $\cK$, then $\varphi(Y_1) = 0$.

\item[{\rm (b)}] $p > r_1$ and $B_1 = (a_1/2)I_{r_1}$.
\end{itemize}
\end{lemma}

\it Proof. \rm
Assume that $A, B$ and $D$ satisfy the assumptions of the lemma.
Then $\|A\| = \|B\| = \|D\| = 1$.
Take $c\in \cS$ such that
$$\|A + B\|_c = \|A + B\| = \|2B + D\| = 2,$$
where $\cS$ is the set in (\ref{c-norm}).
Then
$$2 = \|A + B\|_c
\le \|A\|_c + \|B\|_c \le \|A\| + \|B\| = 2$$
so that $\|A + B\|_c = \|A\|_c + \|B\|_c$.
If $c$ has $p$ nonzero values $c_1\ge \cdots \ge c_p$,
then
$$\sum_{j=1}^p c_j s_j(A+B)
= \sum_{j=1}^p c_j (s_j(A) + s_j(B)).$$
Since
\begin{eqnarray*}
\sum_{j=1}^p c_j s_j(A+B)
&=&\sum_{k=1}^{p-1} (c_k-c_{k+1})\sum_{j=1}^k s_j(A+B)
+ c_p\sum_{j=1}^p s_j(A+B)\\
&\le&
\sum_{k=1}^{p-1} (c_k-c_{k+1})\sum_{j=1}^k [s_j(A)+s_j(B)]
+ c_p\sum_{j=1}^p [s_j(A)+s_j(B)]\\
&=& \sum_{j=1}^p c_j (s_j(A) + s_j(B)),
\end{eqnarray*}
we see that
$\sum_{j=1}^p s_j(A+B) = \sum_{j=1}^p [s_j(A) + s_j(B)]$.
By \cite[Proposition 2.3]{CLT},
there are orthonormal sequences $\{x_1^{(m)}, \ldots, x_p^{(m)}\}$,
$\{y_1^{(m)}, \ldots, y_p^{(m)}\}$ such that
\begin{gather*}\|A\|_c+\|B\|_c = \|A+B\|_c = \lim_m \sum_{j = 1}^p c_j\la (A + B)x_j^{(m)}, y_j^{(m)}\ra\\
\le\ \lim_m \left|\sum_{j = 1}^p c_j\la Ax_j^{(m)}, y_j^{(m)}\ra\right|
+ \lim_m \left|\sum_{j = 1}^p c_j\la Bx_j^{(m)}, y_j^{(m)}\ra\right| \le \|A\|_c + \|B\|_c.\end{gather*}
Thus,
$\lim_m \sum_{j = 1}^p c_j\la Ax_j^{(m)}, y_j^{(m)}\ra = \|A\|_c = \sum_{j = 1}^p c_j s_j(A)$
and similarly for $B$.
Let $w = \min\{p, r_1\}$.
We may assume that for $j = 1, \ldots, w$,
$\la Ax_j^{(m)}, y_j^{(m)}\ra\to s_j(A) = a_1$. Let $P$ be the orthogonal projection of $\cH$
onto $\cK$. As $\cK$ is finite dimensional, we may further assume that $Px_j^{(m)} \to x_j$ and $Py_j^{(m)} \to y_j$ for
$x_j, y_j\in \cK$. Then it is not difficult to show that
$x_j^{(m)}\to x_j$, $y_j^{(m)}\to y_j$ and $x_j = y_j$. Indeed, from
\begin{eqnarray*}
& & \big|\la Ax_j^{(m)}, y_j^{(m)}\ra\big| = \big|\la APx_j^{(m)}, Py_j^{(m)}\ra +
\la A(x_j^{(m)} - Px_j^{(m)}), y_j^{(m)} - Py_j^{(m)}\ra\big|\\
& \le & a_1\|Px_j^{(m)}\| \|Py_j^{(m)}\| + a_2\|x_j^{(m)} - Px_j^{(m)}\| \|y_j^{(m)} - Py_j^{(m)}\|\\
&\le & \big[a_1\|Px_j^{(m)}\|^2 + a_2\|x_j^{(m)} - Px_j^{(m)}\|^2\big]/2
+ \big[a_1\|Py_j^{(m)}\|^2 + a_2\|y_j^{(m)} - Py_j^{(m)}\|^2\big]/2,
\end{eqnarray*}
we get $\|x_j^{(m)} - Px_j^{(m)}\|, \|y_j^{(m)} - Py_j^{(m)}\|\to 0$ and so
$x_j^{(m)}\to x_j$ and $y_j^{(m)}\to y_j$. Then $x_j = y_j$ follows from
$a_1\la x_j, y_j\ra = \la Ax_j, y_j\ra = a_1$.
Moreover, the set $\{x_1, \ldots, x_w\}$
is orthonormal as $\{x_1^{(m)}, \ldots, x_p^{(m)}\}$ is orthonormal
for each $m$. Consider the $p\times p$ matrices $(\la Ax_j^{(m)}, y_k^{(m)}\ra)$ and
$(\la Bx_j^{(m)}, y_k^{(m)}\ra)$.
Passing to a subsequence if needed,
$(\la Ax_j^{(m)}, y_k^{(m)}\ra) \rightarrow \hat A$ and
$(\la Bx_j^{(m)}, y_k^{(m)}\ra) \rightarrow \hat B$ for $p\times p$ matrices $\hat A$ and $\hat B$.
Note that if $p \ge r_1$, the upper left $r_1\times r_1$ principal submatrices
of $\hat A$ and $\hat B$ are $(\la Ax_j, x_k\ra)$ and $(\la Bx_j, x_k\ra)$ respectively.
They are the compressions of $A$ and $B$ onto $\cK$ relative to the orthonormal basis
$\{x_1, \ldots, x_{r_1}\}$. If $p < r_1$, $\hat A$ and $\hat B$
are the compressions of $A$ and $B$ onto the $p$-dimensional subspace of $\cK$
spanned by $\{x_1, \ldots, x_p\}$.
Let $C = \diag(c_1, \dots, c_p)$. Then
$$1 = \|A\|_c = \tr C\hat A\quad \hbox{ and }
\quad 1 = \|B\|_c  = \tr C\hat B.$$

\medskip
We have the following situations.
\begin{itemize}
\item[{\rm (a)}] $p < r_1$. Let $C_1 = C\oplus 0_{r_1 - p}$ and $\varphi$ be
the functional on $\BK$ given by $X\mapsto \tr C_1X$. Then $\varphi(B_1) = \tr C_1B_1 = 1$.
If $Y\in\BH$ satisfies the said condition, then
$$|1 \pm \tr C_1\hat Y_1| = |\tr  C_1(\hat B_1 \pm \hat Y_1)| \le \|\hat B_1 \pm \hat Y_1\|_c \le
\|B \pm Y\|_c\le \|B\pm Y\| = 1.$$
So, $\varphi(Y_1) = \tr(C_1\hat Y_1) = 0$.

\medskip
\item[{\rm (b)}] $p \ge r_1$. Let $\hat D = \hat A - \hat B$.
Then
$$\tr C\hat D =  \tr C(\hat A - \hat B)  = \tr C\hat A
- \tr  C\hat B =\|A\|_c - \|B\|_c = 0.$$
It follows that
$$1 = \tr C(\hat B-\hat D) \le \|\hat B-\hat D\|_c\le \|B - D\|_c \le \|B-D\| = 1.$$
By Lemma \ref{Assertion}, $\hat B - \hat D$ is positive semidefinite.
As the upper left
$r_1\times r_1$ principal submatrix of $\hat B - \hat D$
is the compression of $B - D$ on $\cK$, $B_1 - D_1\ge 0$, where $D_1$ denotes the compression
of $D$ onto $\cK$.
Now, apply the above argument to $D$
starting with $\|A+D\| = 2$, we get a $d\in \cS$
such that $\|A+C\|_d = \|A\|_d + \|C\|_d$.
Suppose $d$ has $q$ nonzero entries.
If $d < r_1$, the compression of $D - B$
onto a $d$ dimensional subspace of $\cK$ is positive
semidefinite with norm one. There is a unit vector $v \in \cK$ such that
$\la Cv, v \ra > \la Bv, v\ra$. This contradicts
our observation $B_1 - C_1\ge 0$.
So, $d\ge r_1$ and we have $C_1 - B_1\ge 0$. It follows that
$B_1 = C_1 = (a_1/2)I_{r_1}$.

Note that $p = r_1$ cannot happen.
Otherwise we will get
$$\|A\| = \sum_{j = 1}^p c_ja_1
= 1\quad\mbox{and}\quad \|B\| = \sum_{j = 1}^p c_j(a_1/2) = 1,$$
which is absurd.
\qed
\end{itemize}

Conditions (a) and (b) in the above lemma can indeed happen as
the following examples show.

\begin{example}
Let $\cH = \ell_2$ and $\BH$ under two different unitarily invariant norms
defined below.
\begin{enumerate}
\item Consider the norm on $\BH$ defined by
$$\|A\| = \max\left\{5s_1(A)/2, s_1(A) + s_2(A) + s_3(A)\right\}.$$
Let $A = (2/5)I_2\oplus (1/5)I$.
Then it is not hard to check that $A$ is an extreme point of the unit ball $\cB$.
We have $A = B+D$ with
$B = (2/5){\rm diag}(1, 0)\oplus (1/10)I$ and
$D = (2/5){\rm diag}(0, 1)\oplus (1/10)I$. Moreover,
$\mu B + \nu D = (2/5){\rm diag}(\mu, \nu)\oplus [(\mu + \nu)/10]I$.
A direct computation shows that
$$\|\mu B + \nu D\| = \max\{|\mu|, |\nu|\}.$$
We have $\|A + B\|_{\rm sp} = \|A + B\| = 2$. So, $(r_1, p) = (2,1)= 2$
in this example.

\medskip
\item Consider $\BH$ under the Ky Fan 2-norm $\|A\|_2 = s_1(A) + s_2(A)$.
Then $A = {\rm diag}(1, 0, \ldots)$ is an extreme point of $\cB$
by \cite[Theorem 3.2]{CLT}.
Then $A = B + D$ with
$B = {\rm diag}(1/2, 1/2, 0, \ldots)$ and $D = {\rm diag}(1/2, -1/2, 0, \ldots)$.
Moreover,
$$\|\mu B + \nu D\|_2 = (1/2)(|\mu + \nu| + |\mu - \nu|) = \max\{|\mu|, |\nu|\}$$
for all
$\mu, \nu\in \mathbb R$.
We have $\|A + B\|_2 = 2$ and so
$(r_1,p) = (1,2)$.\qed
\end{enumerate}
\end{example}

\medskip
Suppose $U$ is a maximal partial isometry. Then $A = U/\|U\|$
is an extreme point of $\cB$. It can be written as $A = B + D$ for which (\ref{norm})
is satisfied. However, the example above shows that it is not a characteristic of
this type of extreme points. Indeed, for every positive integer $N$, we can write $A = A_1 + \cdots + A_N$
such that
\begin{equation}
\label{norm1}\big\|\mbox{$\sum_{j=1}^N a_j A_j$}\big\| = \max\{ |a_j|: 1 \le j \le N\}.\end{equation}
For other extreme points of $\cB$, the number of operators in such a decomposition
is limited by the multiplicity of its largest singular value.

\medskip
\begin{theorem} \label{3.2}
Let $A$ be an extreme point of the norm ball $\cB$.
Then $A$ is a scalar multiple of a maximal partial isometry if and only if
$A = A_1 + \cdots + A_N$ for $N = n^2+1$ such that (\ref{norm1}) is satisfied.
\end{theorem}

\it Proof. \rm For the forward implication, assume that
$\|I\| = f(1, \ldots, 1) = 1$. Let $A$ be a maximal partial isometry.
Then $\|A\| = 1$.
Take orthogonal projections $P_1, \dots, P_N$ each having rank at
least $n$ such that
$P_1 + \cdots + P_N = I$ and put $A_j = AP_j$ or $A_j = P_jA$ for $j =
1, \dots, N$, depending on $A^*A = I$ or $AA^* = I$.
Then $A = \sum_{j=1}^N A_j$. Moreover,
for any $a_1, \ldots, a_N\in \mathbb R$, $\sum_{j=1}^N  a_j A_j$ has the $n$ largest singular values
equal to $|a_\ell|, \dots, |a_\ell|$ if $|a_\ell|= \max\{|a_j|: j = 1,
\dots, N\}$.
Thus, $\|\sum_{j=1}^N a_j A_j\| = |a_\ell| f(1,\dots, 1)= f(I) = |a_\ell|$.

\medskip
For the converse, suppose $A$ is an extreme point of $\cB$
but not a scalar multiple of a maximal partial isometry.
By Lemma \ref{Extreme}, $A = \sum_{j = 1}^n s_j(A) x_jy_j^* + s_n(A)U$,
for $s_j(A), x_j, y_j$'s and $U$ as in the lemma.
Let $V_1$ be given by $x_j\mapsto y_j$ for $j = 1, \ldots, n$ and $V_2$
any unitary operator from $\{x_1, \ldots, x_n\}^\perp$ onto $\{y_1, \ldots, y_n\}^\perp$.
Then $V = V_1\oplus V_2$ is unitary and $VA$ has the form in Lemma \ref{Claim}.

\medskip
Suppose $A = A_1 + \cdots + A_N$ satisfies the said property.
As $V$ is unitary, the singular values of $B$ and $VB$ are the same
for any $B\in \BH$.
So, $VA = VA_1 + \cdots + VA_N$ also satisfies the same property.
For simplicity, let us assume that $A$ is already in that nice form.
Let $B_1, \dots, B_N$ be the compression of $A_j$
onto $\cK = \ker(A - a_1I)$. Note that by Lemma \ref{Extreme}, $\dim \cK < n$.
Now, $A = B + D$ for $B = A_j$ and $D = (A - A_j)$, and (\ref{norm}) holds.
We must have either (a) or (b) of Lemma \ref{Claim}.

\medskip
{\bf Case 1.}
Suppose for each $j$, there is a $\varphi_j\in \BK^*$ such that
$\varphi(B_j) = 1$ and that $\varphi(Y_1) = 0$ whenever
$Y\in \BH$ satisfies
$\|\mu A_j + \nu Y\| = \max\{|\mu|, |\nu|\}$
for any $\mu, \nu \in \IR$ and $Y_1$ is the compression of $Y$ onto $\cK$.
The condition is clearly satisfied by $Y = A_k$ for $k\ne j$.
So, the bi-orthogonal condition $\varphi_j(B_k) = \delta_{jk}$ ensures
$B_1, \dots, B_N$ are linearly independent. However, $\dim \cK \le n - 1$.
There cannot be so many linearly independent matrices.
This case is impossible.

\medskip
{\bf Case 2.} Suppose for some $A_j$, say $A_1$,
there is a $c\in \cS$ with nonzero entries $c_1, \ldots, c_p$
for $p\ge r_1$ such that $\|A + A_1\|_c = \|A + A_1\| = 2$.
By Lemma \ref{Claim}, the compression of $A_1$ to $\cK$ is $(a_1/2)I_{r_1}$.

\medskip
For $j > 1$, we construct the $p\times p$ matrices $\hat A$, $\hat A_1$ and
$\hat A_j$ as in the proof of Lemma \ref{Claim}.  Then for $C = {\rm diag}(c_1, \ldots, c_p)$,
$$|1 \pm \tr C\hat A_j| = |\tr  C(\hat A_1 \pm \hat A_j)| \le \|\hat A_1 \pm \hat A_j\|_c \le
\|A_1 \pm A_j\| = 1.$$
Thus $\tr(C\hat A_j) = 0$. It follows that
$$2 =  \tr C[A + (\hat A_1 + \hat A_j)] \le \|\hat A + (\hat A_1 + \hat A_j)\|_c
\le \|A + (A_1 + A_j)\|_c \le \|A + (A_1 + A_j)\| = 2.$$
We have $\|A + (A_1 + A_j)\|_c = \|A + (A_1 + A_j)\| = 2$ for $c$ with $p\ge r_1$
nonzero entries.
So, Lemma \ref{Claim} applied to
$B = A_1 + A_j$ and $D = A - B$ gives us $B_1 + B_j = (a_1/2)I_{r_1}$.
(Here, note that $\dim \cK < n$ ensures $n > 1$ so that $n^2+1 \ge 5$
and $A-B$ is nontrivial.)
But we already have
$B_1 = (a_1/2)I_{r_1}$.
So, $B_j = 0_{r_1}$ for all $j > 1$.
Then the compression of $A$ onto $\cK$ is
$B_1 + \cdots + B_N = (a_1/2)I_{r_1}$,
which contradicts the form of $A$. This case is also impossible.
\qed

\medskip
\it Proof of Theorem \ref{isometry}.
\rm
The sufficiency is clear. To prove the necessity, suppose $L: \BH\to \BH$ satisfies the
hypothesis of the theorem.
Observe that by the Mazur-Ulam Theorem, the map $\hat L: A \mapsto L(A)-L(0)$ is
real linear and satisfies
$\|\hat L(A)\| = \|A\|$ for all $A \in \BH$.
So, $\hat L$ preserves the extreme points of $\cB$.
By Theorem \ref{3.2},
$\hat L$ sends maximal isometries to maximal isometries.
It follows from \cite[Solution to Problem 107]{H} (see also \cite[Lemma 3]{R}) that $\hat L$ is indeed
an isometry of the spectral norm.
Using the characterization \cite[Corollary 3.3]{Dang}, we conclude that $L$ is of the asserted form.
\qed

\medskip
\section*{Acknowledgement}

Li is an affiliate member of the Institute for Quantum
Computing, University of Waterloo; his
research was partially supported by the
Simons Foundation Grant 851334. The authors would like to
thank Professor Ngai-Ching Wong for some inspiring discussions and comments,
and also the anonymous referee for the helpful comments.

\end{document}